\documentclass[10.9pt,a4paper]{amsart}
\usepackage{a4wide}

\usepackage[utf8]{inputenc}
\usepackage[english]{babel}

\usepackage{amsfonts,amssymb,amsmath,pinlabel,array,hhline}
\usepackage[nobysame,alphabetic, initials]{amsrefs}

\usepackage[dvipsnames]{xcolor}
\usepackage{xcolor}
\usepackage[colorlinks,
    linkcolor={red!50!black},
    citecolor={blue!50!black},
    urlcolor={blue!80!black}]{hyperref}

\usepackage[position=b]{subcaption}

\usepackage{tikz}
\usepackage[all]{xy}
\usepackage{enumitem}
\makeatletter
\newcommand{\mylabel}[2]{#2\def\@currentlabel{#2}\label{#1}}
\usetikzlibrary{calc, matrix, arrows, cd}

\newcommand{\bsm}{\left(\begin{smallmatrix}}
\newcommand{\esm}{\end{smallmatrix}\right)}

\newtheorem{theorem}{Theorem}[section]

\newtheorem{addendum}[theorem]{Addendum}

\newtheorem{meta}[theorem]{Meta-Statement}
\newtheorem{proposition}[theorem]{Proposition}

\theoremstyle{definition}
\newtheorem{definition}[theorem]{Definition}

\newtheorem{example}[theorem]{Example}

\newtheorem{remark}[theorem]{Remark}
\newtheorem{construction}[theorem]{Construction}
\newtheorem*{claim}{Claim}
\newtheorem*{claim*}{Claim}

\newcommand{\Z}{\mathbb{Z}}

\newcommand{\C}{\mathbb{C}}

\newcommand{\Hom}{\operatorname{Hom}}
\newcommand{\ks}{\operatorname{ks}}
\newcommand{\id}{\operatorname{id}}

\begin{document}
\title{On classifying knotted surfaces}
\author{Anthony Conway}
\address{The University of Texas at Austin, Austin TX 78712}
%, United States}
\email{anthony.conway@austin.utexas.edu}
\begin{abstract}
These notes, which are based on four lectures delivered at the summer school ``Knotted Surfaces in Four-Manifolds" at PCMI in~2026, discuss classifications of locally flat surfaces in simply-connected $4$-manifolds.
%, with a focus on the case where the complements of the spheres have abelian fundamental group.
\end{abstract}

\def\subjclassname{\textup{2020} Mathematics Subject Classification}
\expandafter\let\csname subjclassname@1991\endcsname=\subjclassname
\maketitle

The aim of these lectures is to survey known classifications of locally flat surfaces in simply-connected~$4$-manifolds.
This area can be thought of as lying at the intersection of the study of topological $4$-manifolds with the theory of knotted surfaces.
In essence, in contrast with the smooth category,  the availability of tools from high dimensional topology makes it possible, in principle, to obtain broad structural results for such knotted surfaces.

\medbreak

We briefly outline the history 
%and evolution 
of the subject. 
The topic took off with the work of Freedman and Quinn~\cite{Freedman,QuinnEndsOfMapsIII,FreedmanQuinn} thanks to which it became known that Alexander polynomial one knots are slice and that if an embedded~$2$-sphere~$S \subset S^4$ satisfies~$\pi_1(S^4 \setminus S)=\Z$, then it is unknotted.
Shortly afterwards,  Lawson proved that if a projective plane~$S \subset S^4$ satisfies~$\pi_1(S^4 \setminus S)=\Z_2$, then it is unknotted (see also~\cite{FinashinKreckViro,KreckOnTheHomeomorphism}) and, in the early nineties, Lee and Wilczy\'nski initiated a thorough study of closed orientable surfaces whose complement have finite cyclic fundamental group~\cite{LeeWilczyOdd,LeeWilczy,LeeWilczyGenus} (see also~\cite{HambletonKreck2}), whereas Boyer classified closed orientable surfaces whose complement is simply-connected~\cite{BoyerRealization}.
During the aughts,  attention shifted to classical knot concordance, where high dimensional techniques have less sway (nonetheless see~\cite{FriedlTeichner,Davis} for striking applications of surgery theory to this area).
%but, 
Over the recent years, there has been a resurgence in attempting to classify knotted surfaces~\cite{ConwayPowellDiscs,FellerLewarkBalanced,FellerMillerNagelOrsonPowellRay, ConwayPowell,OrsonPowellSpines,  KasprowskiPowellRayTeichnerEmbedding, ConwayPiccirilloPowell,ConwayDiscs, ConwayOrsonPowell, ConwayMiller, ConwayOrson, ConwayOrsonPencovitch,ConwayGalvin, Pyronneau,ConwayKasprowskiSurfaces, ConwayUnknotting, PencovitchSlice}. 
%KasprowskiPowellRayGluck It's not really a surface result.

\medbreak

The organisation of this lecture series is as follows.
Lecture~\ref{sec:1} makes precise what we mean by a ``classification" and reviews the aforementioned early successes of Freedman,  Lawson, Boyer and Lee-Wilczy\'nski. 
Lecture~\ref{sec:2} collects some of the flagship theorems from $4$-manifold topology such as the $s$-cobordism theorem and Freedman's classification of closed simply-connected $4$-manifolds.
Lecture~\ref{sec:3} focuses on uniqueness results for knotted surfaces. Finally, Lecture~\ref{sec:4} concludes with realisation results.
Some exercises are listed at the end of the notes.

\medbreak

Since these lectures took place during the second week of a three week long summer school devoted to knotted surfaces,  some familiarity is assumed with the basics of knotted surface theory.
Some knowledge of algebraic topology is also assumed.
Finally,  we do not aim to survey all that is known about locally flat surfaces and their invariants, instead we focus on classification results.

\subsection*{Conventions}
We will work exclusively in the topological category with locally flat embeddings.
Manifolds are assumed to be compact and connected and $4$-manifolds are assumed to be oriented.

\subsection*{Acknowledgments}

I wish to thank the organisers of the 2026 PCMI Graduate Summer School on ``Knotted Surfaces in Four-Manifolds" for the opportunity to deliver these lectures.
This work was partially supported by the NSF grant DMS~2303674.

\subsection*{AI Statement}
No AI was used during the preparation of this work.

\newpage

\section{Lecture one: What is a classification?}
\label{sec:1}

In dimension $4$,  in the topological category, tools from high-dimensional topology are available when the fundamental group of the $4$-manifolds is ``good"~\cite{Freedman,QuinnEndsOfMapsIII,FreedmanQuinn}.
In practice, this means that classification results become, in principle, achievable.
The goal of this first lecture is to make these expectations more precise as well as to introduce some first surface invariants.
\medbreak
To begin with,  we agree on what it means for two knotted surfaces to be ``the same".

\begin{definition}
\label{def:Isotopy}
Two embedded surfaces~$S_0 \subset X$ and~$S_1 \subset X$ are 
 \emph{equivalent} if there is a homeomorphism~$\Phi \colon X \to X$ with~$\Phi(S_0)=S_1$ and \emph{isotopic} if, moreover,~$\Phi$ is isotopic to the identity.
\end{definition}

Our next task is to elucidate what is meant by a classification.
Traditionally,  this is understood to mean that there exists a complete set of algebraic invariants that can tell the given topological objects apart.
For surfaces, we phrase this as follows.

\begin{meta}[Uniqueness]
\label{meta:Uniqueness}
Let $X$ be a $4$-manifold.
For every surface $S \subset~X$, there is a collection of algebraic invariants $I(S)$ taking values in a set $\mathcal{A}$ such that if two surfaces~$S_0,S_1 \subset X$ share the same invariants, i.e.~$I(S_0)=I(S_1)$, then~$S_0,S_1$ are equivalent or isotopic.
%, possibly rel. boundary.
\end{meta}

Here, the word ``algebraic" is included in an attempt to avoid tautological statements involving the surface exterior or the surface itself.
We emphasise that both the $4$-manifold and the surfaces may or may not have boundary.
In the case with boundary, we will assume the surfaces are properly embedded so that $\partial S \subset \partial X$.
We do consider the boundary of a surface to be an invariant.

While the statement above reflects the most common understanding of the word ``classification",  a second question is just as important: namely whether or not the algebraic objects from the classification necessarily arise as surface invariants.

\begin{meta}[Existence]
\label{meta:Existence}
Let $X$ be a $4$-manifold and let $\mathcal{A}$ be a complete set of algebraic invariants as in Statement~\ref{meta:Uniqueness}.
For every~$a \in \mathcal{A}$, there exists a surface~$S \subset X$ with~$I(S)=a$.
\end{meta}

The combination of these two statements yields a bijection between the set of equivalence (or isotopy) classes of surfaces in~$X$ and the purely algebraic set $\mathcal{A}$.
In practice, it is common to fix a couple of first invariants, such as the (possibly nonorientable) genus, the boundary,  and the \emph{knot group}~$\pi_1(X\setminus S)$ and to then ask for classifications of surfaces with the given invariants.
For brevity,  it is equally common to refer to surfaces with knot group $\pi$ as \emph{$\pi$-surfaces.}

\begin{remark}
These lectures focus on classifications up to equivalence and isotopy, but other equivalence relations are also the subject of much research.
A prime example is the classification of knotted surfaces up to concordance, see e.g.~\cite{Kervaire,BlanloeilSaeki,Sunukjian,MillerConcordanceAnalogue,
KlugMillerMoreGeneral,DaherPowell,Hellsten}.
\end{remark}

\begin{remark}
We briefly comment on another invariant that might also have come to the reader's mind.
The \emph{homology class} of a closed oriented surface $S \subset X$ refers to the class~$[S] \in H_2(X)$.\
When the surfaces have boundary (resp. are nonorientable) relative homology (resp. homology with~$\Z_2$ coefficients) is to be used instead.
The homology class is an isotopy invariant but not necessarily an equivalence invariant. 
Equivalence invariants can nevertheless be extracted from the homology class of $S$.
\begin{itemize}
\item  When $S$ is a closed orientable surface and~$X$ is assumed to be simply-connected for simplicity, the \emph{divisibility} of~$[S]$, i.e. the largest integer~$d \geq 0$ dividing~$[S] \in H_2(X) \cong \Z^{b_2(X)}$, is invariant under equivalences.
Writing $X_S:=X \setminus \nu S$ for the complement of a tubular neighborhood of $S$, a suggested exercise is to show that  a divisibility~$d$ surface~$S\subset X$ satisfies~$H_1(X_S) \cong \Z_d$.
%and thus, if~$S$ is simple,  that~$\pi_1(X_S)\cong \Z_d$.
For nonorientable surfaces,  the divisibility is defined to be~$d=2$ if~$[S] \in H_2(X;\Z_2)$ is trivial and~$d=1$ if~$[S] \in H_2(X;\Z_2)$ is nontrivial, so that the isomorphism~$H_1(X_S) \cong \Z_d$ still holds.
Divisibility~$d=1$ classes are called \emph{primitive}.
\item A homology class $w \in H_2(X)$ is \emph{characteristic} if $w.x=x.x$ mod~$2$ for every $x \in H_2(X)$ and \emph{ordinary} otherwise.
A closed orientable surface is \emph{characteristic} (resp. ordinary) if $[S]$ is characteristic (resp. ordinary).
Being characteristic/ordinary is once again an equivalence invariant.
%S_1. x=F(S_0).x=S_0 \cdot y=y.y=x.x for every x.
Suitable generalisations exists for surfaces with nonempty boundary and for nonorientable surfaces; the details are left to the reader.
\item The \emph{Euler number} of a closed orientable surface $S\subset X$ is
$$ e(S)=[S] \cdot [S] \in \Z.$$
This integer can also be defined by counting the intersections between $S$ and a generic push-off of~$S$ in the total space of the normal bundle $\nu(S) \subset X$.
This latter definition generalises to nonorientable surfaces; see e.g.~\cite[Section 2]{ConwayOrsonPowell} for a detailed discussion.
\end{itemize}
\end{remark}

Next, we describe some early successes in classifying knotted surfaces,  with a focus on results obtained during the 80s and 90s. 
In order to appreciate the following theorem, it is helpful to recall that unknotted orientable (resp. nonorientable) surfaces in~$S^4$ have knot group~$\Z$ (resp.~$\Z_2$),  as can be seen for example by drawing handle diagrams of the surface exteriors.

\begin{theorem}[{Freedman~\cite{FreedmanICM}}]
\label{thm:FreedmanUnknotting}
Any $\Z$-sphere in $S^4$ is unknotted.
\end{theorem}

Freedman's result therefore shows that a sphere $S \subset S^4$ is unknotted if and only if its knot group is infinite cyclic.
In the framework of our meta-statements,  the theorem states that the set of isotopy classes of~$\Z$-spheres in~$S^4$ is a singleton consisting of the unknotted~$2$-sphere.
Shortly after Freedman's unknotting theorem,  Lawson proved the first ``nonorientable unknotting theorem".

\begin{theorem}[{Lawson~\cite{Lawson}}]
\label{thm:LawsonUnknotting}
Any $\Z_2$-projective plane in $S^4$ is unknotted.
\end{theorem}

 Lawson's result implies that a projective plane $S \subset S^4$ is unknotted if and only if its knot group has order $2$.
In particular, up to isotopy,  there are two $\Z_2$-projective projective planes in~$S^4$, distinguished by their Euler number.
Indeed, Massey proved that the Euler number of a closed nonorientable genus $h$ surface in $S^4$ lies in the range $\{-2h,-2h+4,\ldots,2h-4,2h\}$, see~\cite{Massey}.
Note that each of these values is realised by an unknotted surface, i.e. by a connected sum of unknotted projective planes. 

The next result,  due to Boyer,  provides the classification of surfaces with trivial knot group.

\begin{theorem}[{Boyer~\cite{BoyerRealization}}]
\label{thm:Boyer}
Fix a closed simply-connected~$4$-manifold~$X$ and an integer~$g \geq 0$.
\begin{enumerate}
\item Uniqueness: Genus~$g$ surfaces in~$X$ with simply-connected complements are isotopic if and only if they are homologous.
\item Realisation: A homology class $x \in H_2(X)$ is realised by a genus~$g$ surface with simply-connected complement if and only if it is primitive and, when~$g=0$ and $x$ is characteristic, it additionally satisfies the \emph{Kervaire-Milnor condition}:
\begin{equation}
\label{eq:KM}
\ks(X)=\frac{1}{8}(\sigma(X)-x\cdot x) \in \Z_2.   \tag{km}
\end{equation}
\end{enumerate}
\end{theorem}

Here~$\sigma(X) \in \Z$ denotes the signature of~$X$ and $\ks(X) \in \Z_2$ stands for the Kirby-Siebenmann invariant of $X$.
The latter can be defined as
$$
\ks(X)=
\begin{cases}
0 \quad &\text{if~$X \#_n S^2 \times S^2$ is smoothable for some~$n \geq 0$} \\
1  \quad &\text{otherwise}.
\end{cases}
$$
The reason we refer to~\eqref{eq:KM} as the \emph{Kervaire-Milnor condition} is as follows.
Kervaire and Milnor proved that if a homology class~$x \in H_2(X)$ can be represented by a smoothly embedded sphere, then~\eqref{eq:KM} holds~\cite{KervaireMilnor}; the corresponding result in the topological category is due to Freedman-Kirby~\cite{FreedmanKirby}.
For example, the class~$3 \in H_2(\C P^2) \cong \Z$ cannot be represented by a sphere.

Theorem~\ref{thm:Boyer} therefore implies that genus~$g$ surfaces in~$X$ with trivial knot group correspond bijectively with primitive homology classes in~$X$ that additionally satisfy the Kervaire-Milnor condition when~$g=0$ and the class is characteristic.

\medbreak

Another striking result from the nineties is Lee-Wilczy\'nski's and Hambleton-Kreck's classification of knotted spheres with finite cyclic knot group.
We first state the outcome of their work and then introduce the definitions needed to fully make sense of it.

\begin{theorem}[{Lee-Wilczy\'nski~\cite{LeeWilczyOdd,LeeWilczy,LeeWilczyGenus} and Hambleton-Kreck~\cite{HambletonKreck2}}]
\label{thm:LW}
Fix a closed simply-connected~$4$-manifold~$X$ and a positive integer $d>0$.
\begin{enumerate}
\item  
Two~$\Z_d$-spheres in~$X$ are equivalent
%\footnote{In fact, under certain assumptions on $X$,  Lee and Wilczynski show that~$\Z_d$-spheres in~$X$ are isotopic if and only if they are homologous~\cite[Theorem 1.1]{LeeWilczy}.} 
if and only if they are either both ordinary or both characteristic,  they have the same Euler number, and the pointed intersection forms of their branched covers are isometric.
Given a homology class $x \in H_2(X)$ of divisibility $d$, when~$d=1$ or when~$X$ satisfies both~$b_2(X)>6$\footnote{
Under additional hypotheses, such results are also available when~$3 \leq b_2(X) \leq 6$; see~\cite[Addendum 1]{LeeWilczy} as well as~\cite[Theorem~4.5]{HambletonKreck2}.
In both references,~\eqref{eq:UniquenessIneq} is still required to hold.
}
 and
\begin{equation}
\label{eq:UniquenessIneq}
 b_2(X) > \underset{0 \leq j <d}{\operatorname{max}} \ \Big| \sigma(X)-\frac{2j(d-j)}{d^2}x\cdot x \Big|, 
 \tag{\textgreater}
 \end{equation}
two $\Z_d$-spheres representing $x$ are isotopic.
% if and only if they are homologous.
\item  A homology class $x \in H_2(X)$ is realised by a~$\Z_d$-surface of genus~$g \geq 0$ if and only if its divisibility equals~$d$, it satisfies the \emph{Rochlin inequality}
\begin{equation}
\label{eq:Rochlin}
b_2(X) +2g \geq 
\underset{0 \leq j <d}{\operatorname{max}}
\Big|
\sigma(X) -\frac{2j(d-j)}{d^2} x\cdot x 
\Big|,   \tag{Ro}
\end{equation}
and, when~$g=0$ and $x$ is characteristic, it additionally satisfies the Kervaire-Milnor condition from~\eqref{eq:KM}.
\end{enumerate}
\end{theorem}

Before providing details on equivariant intersection forms, we note that we refer to~\eqref{eq:Rochlin} as the \emph{Rochlin inequality} because Rochlin proved that if a homology class $x \in H_2(X)$ can be represented by a smoothly embedded orientable genus $g$ surface, then the inequality~\eqref{eq:Rochlin} holds when $d$ is a prime power~\cite{Rochlin}.
%%Without any pi_1-assumption.
The fact that Rochlin's result also holds in the topological category was known to researchers working in the area after Freedman and Quinn's breakthroughs.

%%%Could add a condition on the isometry, technically. Bof.
%\begin{remark}
%We note that if $X$ satisfies~\eqref{eq:UniquenessIneq}, then it cannot be definite. 
%In particular,  up to isotopy, complete classifications of surfaces with cyclic knot group remain unknown in general, even in the case of spheres.
%\end{remark}

Next, we review equivariant intersection forms and, in particular,  the invariant that appears in the first statement of Lee and Wilczy\'nski's Theorem~\ref{thm:LW}.

\begin{construction}
\label{cons:Forms}
The \emph{equivariant intersection form} of a~$4$-manifold~$W$ with fundamental group~$\pi:=\pi_1(W)$ refers to the hermitian form
\begin{align}
\label{eq:Equivariant}
\lambda_{\widetilde{W}} \colon H_2(\widetilde{W}) \times H_2(\widetilde{W}) &\to \Z[\pi]  \nonumber \\
(x,y)& \mapsto \sum_{g \in \pi} (x \cdot gy) g^{-1}. 
\end{align}
Here, ~$\widetilde{W}$ denotes the universal cover of~$W$ and being \emph{hermitian} means that~$\lambda_{\widetilde{W}}(y,x)=\overline{\lambda_{\widetilde{W}}(x,y)}$, where the overline refers to replacing all instances of~$g \in \pi$ by~$g^{-1}$. 
The sum in~\eqref{eq:Equivariant} is finite because homology classes are compactly supported.
We note that the equivariant intersection form also admits a homological definition in terms of Poincar\'e duality, but we will not dwell on this fact in these notes.
Typically, we will be considering the case where $W=X\setminus \nu(S)=:X_S$ is the complement of a tubular neighborhood of a surface~$S \subset X$.
In this case,  the (isometry type of the) equivariant intersection form is an equivalence-invariant of the surface.
%We note that the existence of tubular neighborhoods in this dimension is a deep fact~\cite[Theorem 9.3]{FreedmanQuinn}.

The same definition as above can be made for branched covers,  but we spell it out explicitly given its relevance to Theorem~\ref{thm:LWRealisation}.
Write $\Z_d=\langle T \mid T^d=1\rangle$.
For a surface~$S \subset X$ with~$H_1(X_S) \cong \Z_d$, the equivariant intersection form of the $d$-fold branched cover $\Sigma_d(S)$ refers to the hermitian form
\begin{align*}
\lambda_{\Sigma_d(S)} \colon H_2(\Sigma_d(S)) \times H_2(\Sigma_d(S)) &\to \Z[\Z_d] \\
(x,y) &\mapsto \sum_{k=0}^{d-1} Q_{\Sigma_d(S)}(x,T^ky) T^{-k}.
\end{align*}
The \emph{pointed equivariant intersection form} of $S$ refers to the pair~$(\lambda_{\Sigma_d(S)},[\widetilde{S}])$, where~$\widetilde{S} \subset \Sigma_d(S)$ denotes the branched set of the cover~$\Sigma_d(S) \to X$.
An \emph{isometry} of such pointed forms is an isometry of the underlying hermitian forms that takes one homology class to the other.
\end{construction}

\begin{example}
\label{ex:Equivariant}
Explicit calculations of equivariant intersection forms typically require handle decompositions of the surface exteriors.
We will not go through such calculations explicitly, but list some of the outcomes to build intuition.
\begin{enumerate}
\item The equivariant intersection form of the exterior of the unknotted surface of genus~$g$, which we denote~$U \subset S^4$, is given by
$$
(H_2(\widetilde{X}_U),\lambda_{\widetilde{X}_U})
\cong
\left(
\Z[t^{\pm 1}]^{\oplus 2g},
\bsm
0&t-1 \\
t^{-1}-1&0
\esm^{\oplus g}
\right).
$$
This fact can be deduced by considering a handle diagram for the exterior of the unknotted surface and passing to the infinite cyclic cover~\cite[Lemma 6.1]{ConwayPowell}. 
In fact, 
%swhen $g \neq 1,2$,  
arduous algebra shows that the equivariant intersection form of any $\Z$-surface~$S \subset S^4$ is isometric to~$\bsm
0&t-1 \\
t^{-1}-1&0
\esm^{\oplus g}$; see~\cite[Theorem 7.4]{ConwayPowell} and~\cite[Theorems 1.6 and 1.7]{ConwayUnknotting}.
\item For any~$\Z_2$-sphere~$S \subset \C P^2$,  it is a suggested exercise to show that
$$
(H_2(\widetilde{X}_S),Q_{\widetilde{X}_S})
\cong 
(\Z_-,\pm 2).
$$
Here,~$\Z_-$ denotes the~$\Z[\Z_2]$-module whose underlying abelian group is~$\Z$ and where the action is induced by~$Tx=-x$ for every~$x \in \Z$.
The long exact sequence of the pair then implies $H_1(\partial \widetilde{X}_S)\cong \Z_2$ and it also follows that the equivariant intersection form can be described as
$$
(H_2(\widetilde{X}_S),\lambda_{\widetilde{X}_S})
\cong 
(\Z_-,(x,y) \mapsto \pm 2(1-T)xy).
$$
Some additional care shows that the sign is determined by the Euler number and therefore by the homology class of $S$.
Describing the \emph{pointed} intersection form in even this simple case is already arduous, especially without using handle decompositions.
\end{enumerate}
\end{example}

We conclude this section with a further discussion of Rochlin's theorem.
\begin{remark}
\label{rem:Rochlin}
Given a closed orientable characteristic surface $S \subset X$ and a loop~$\gamma \subset S$,  \emph{the Rochlin quadratic form}
$$q_S \colon H_1(S;\Z_2) \to \Z_2$$
assigns to a class~$[\gamma] \in H_1(S;\Z_2)$ the mod~$2$ integer~$S \cdot D +W(D)$, where~$D \subset X$ is an embedded disc with~$\partial D=\gamma=D \cap S$ and~$W(D) \in \Z$ is the obstruction to the surface-induced framing on~$\gamma$ extending over the disc~$D$.
%the disc framing on~$\gamma$ agreeing with the surface-induced framing on it.
%a framing obstruction.
We refer to e.g.~\cite{FreedmanKirby,Matsumoto} for the proof that~$q_S$ is well defined and that it provides a quadratic refinement of the intersection form on $H_1(S;\Z_2)$.
In this setting,  Freedman and Kirby~\cite{FreedmanKirby} proved that 
$$ 
\operatorname{Arf}(q_S)
\equiv \ks(X) -\frac{1}{8}(\sigma(X)-[S] \cdot [S]) 
\mod 2.
$$
Here, the \emph{Arf invariant} $\operatorname{Arf}(q_S) \in \Z_2$ can be defined as the value taken most often by $q_S$.
When $X$ is spin,  this congruence shows that the Kirby-Siebenmann invariant is determined by the signature:
\begin{equation}
\label{eq:KSSpin}
\ks(X) \equiv  \frac{1}{8} \sigma(X) \mod 2.
\end{equation}
This can be seen by taking $S\subset X$ to be a small (necessarily nullhomologous) unknotted sphere.

There is an analogue of the Rochlin form in the nonorientable setting due to Guillou-Marin~\cite{GuillouMarin} where, this time, the form $q_S$ is $\Z_4$-valued.
For generalisations of these congruences to surfaces with nonempty boundary, we refer to~\cite{Klug}.
\end{remark}

\section{Lecture two: Tools of the trade}
\label{sec:2}

This lecture focuses on three fantastic theorems that power much of $4$-manifold topology in the topological category.
We conclude with some first examples that illustrate how these results can be applied to study knotted surfaces.

\subsection{The~$s$-cobordism theorem}

Recall that an \emph{$h$-cobordism} between two closed~$n$-manifolds~$M_0$ and~$M_1$ refers to an~$(n+1)$-dimensional cobordism~$(W;M_0,M_1)$ such that the inclusion~$M_i \to~W$ is a homotopy equivalence for~$i=0,1$.
When $n \geq 5$, the (topological)~$h$-cobordism theorem states that if~$M_0$ and~$M_1$ are simply-connected, then~$W$ is homeomorphic rel.~$M_0$ to the product cobordism~$M_0 \times [0,1]$; see~\cite{Newman,KirbySiebenmann}.
%The work of Freedman,  Quinn and Freedman-Quinn~\cite{Freedman,QuinnEndsOfMapsIII,FreedmanQuinn} ensures that the same conclusion holds in dimension~$4$.

The statement in the non simply-connected case is known as \emph{the~$s$-cobordism theorem.}
Without going into the details,  in dimensions $n \geq 5$, it states that the conclusion of the~$h$-cobordism theorem still holds if~$W$ is an~\emph{$s$-cobordism} i.e.  if the inclusions~$M_0 \to W$ and~$M_1 \to W$ are \emph{simple} homotopy equivalences.
We refer to e.g.~\cite{LueckMacko} for an account of simple homotopy equivalences, the Whitehead group and Whitehead torsion in the context of the $s$-cobordism theorem.

The~$s$-cobordism theorem also holds for manifolds with boundary; see e.g.~\cite[page 4]{KirbySiebenmann}.
To make sense of the statement, recall that given a homeomorphism~$f \colon \partial M_0 \to \partial M_1$,   a \emph{cobordism~$(W;M_0,M_1)$ rel.~$f$} refers to~$(n+1)$-manifold~$W$ with~$\partial W \cong M_0 \cup_{\id} (\partial M_0 \times [0,1]) \cup_f M_1$.
Again, such a cobordism is \emph{an~$s$-cobordism rel.~$f$} if the inclusion~$M_i \to W$ is a simple homotopy equivalence for $i=0,1$.

The work of Freedman,  Quinn and Freedman-Quinn~\cite{Freedman,QuinnEndsOfMapsIII,FreedmanQuinn} ensures that the $s$-cobordism theorem still holds in dimension~$4$ provided the fundamental group is well behaved.

\begin{theorem}[The~$s$-cobordism theorem rel. boundary]
\label{thm:scob}
Let $f \colon \partial M_0 \to \partial M_1$ be a homeomorphism.
For~$n \geq 5$, any~$(n+1)$-dimensional~$s$-cobordism~$(W;M_0,M_1)$ is homeomorphic rel.~$M_0 \cup (\partial M_0\times [0,1])$ to the product~$s$-cobordism $(M_0 \times [0,1];M_0,M_1)$.
The same statement holds if~$n=4$ and~$\pi$ is a good group.
\end{theorem}

The definition of a good group involves certain geometric objects known as \emph{gropes}.
We omit the precise definition (referring to~\cite{FreedmanTeichner1} or~\cite[Chapter 19]{DET} for the details) but record that solvable groups are good.
In particular,  finite groups and abelian groups are good.
Finally,  it is a suggested exercise to show that if the conclusion of Theorem~\ref{thm:scob} holds, then $h\colon \partial M_0 \to \partial M_1$ extends to a homeomorphism $M_0 \to M_1$.

\subsection{Surgery theory}

Since, under the assumptions of Theorem~\ref{thm:scob},~$s$-cobordant~$n$-manifolds are homeomorphic, it becomes natural to ask: how can one find an~$s$-cobordism rel.~$f$ between two given~$n$-manifolds?
When the homeomorphism~$f \colon \partial M_0 \to \partial M_1$ extends to a homotopy equivalence~$M_0 \to M_1$, surgery theory suggests a programme for promoting this homotopy equivalence to a homeomorphism. 
A second approach known as \emph{modified surgery theory} consists of finding a particularly well behaved cobordism~$(W;M_0,M_1)$ rel.~$f$ and then surgering it to an~$s$-cobordism.

Classical surgery theory was developed by Browder~\cite{BrowderHomotopyType,BrowderSimply}, Novikov~\cite{NovikovHomotopically}, Sullivan~\cite{SullivanThesis} and Wall~\cite{WallSurgeryOnCompact} in the sixties, whereas modified surgery theory was developed by Kreck in the late seventies~\cite{KreckUnpublished,KreckSurgeryAndDuality}.
Both approaches are technically demanding, and become increasingly unwieldy  the further one veers away from the simply-connected case.
Since the focus of these notes is on knotted surfaces, we will content ourselves with stating one of the outcomes of classical surgery when applied to $4$-manifolds.

\begin{theorem}
\label{thm:SurgeryTheory}
Assume that~$\pi:=\pi_1(X_1)$ is trivial,  finite cyclic, or a solvable Baumslag--Solitar group.
Let~$X_0$ and~$X_1$ be~$4$-manifolds with fundamental group~$\pi$ and~$\pi_1(\partial X_i) \to \pi_1(X_i)=:\pi$ surjective for~$i=0,1$.
If~$\ks(X_0)=\ks(X_1)$ and a homeomorphism~$f \colon \partial X_0 \to \partial X_1$ extends to a  homotopy equivalence~$X_0 \to X_1$,  then~$X_0$ and~$X_1$ are~$s$-cobordant rel.~$f$.
\end{theorem} 

This type of result is known to researchers working in the area (with the simply-connected case dating back to the work of Freedman~\cite{Freedman}); see e.g. the references within~\cite[Proof of Proposition 1.5]{ConwayKasprowski4Manifolds}.
%For finite cyclic, there is SK trickery (Wh nontrivial).
We emphasise that in Theorem~\ref{thm:SurgeryTheory} the restrictions on the fundamental group are not merely stated for convenience: even for relatively straightforward non-cyclic groups,  surgery theory typically becomes challenging to apply. 
The arduous question of extending homeomorphisms to homotopy equivalences will be briefly discussed in the next lecture in the context of surface exteriors; we refer to~\cite{ConwayKasprowski4Manifolds} for the broader theory and for further details.

\medbreak

We recall how Theorem~\ref{thm:SurgeryTheory} combines with Whitehead's homotopy classification and the s-cobordism theorem to yield the homeomorphism classification of closed simply-connected~$4$-manifolds.
 
\begin{theorem}[The classification of closed,  simply-connected~$4$-manifolds~\cite{Freedman}]
\label{thm:HomeoClassification}
~
\begin{itemize}
\item Two closed,  simply-connected~$4$-manifolds are homeomorphic if and only if they have isometric intersection forms and the same Kirby-Siebenmann invariant.
\item For every nonsingular, symmetric,  bilinear form~$Q \colon H \times H \to \Z$, there exists a closed simply-connected~$4$-manifold~$X$ with~$Q_X \cong Q$.
When $Q$ is even, $\ks(X)=\operatorname{sign}(Q)/8$; when~$Q$ is odd, either value of the Kirby-Siebenmann invariant can be realised.
\end{itemize}
\end{theorem}
\begin{proof}
Milnor (building on work of Whitehead~\cite{Whitehead4manifold}) proved that smooth simply-connected~$4$-manifolds are homotopy equivalent if and only if they have isometric intersection forms~\cite{Milnor4Manifolds}.
A proof that also works in the topological category is given
in~\cite[Chapter V, Theorem 1.5]{HusemollerMilnor}.
%The same assertion holds in the topological category because every closed $4$-manifolds is homotopy equivalent to a $4$-dimensional CW complex; see~\cite[Theorem 3.16]{FriedlNagelOrsonPowell} for an exposition of this result.
The first statement then follows from Theorem~\ref{thm:SurgeryTheory} and the $h$-cobordism theorem. 
We outline the proof of the second statement.
Form a framed link~$L \subset S^3$ whose framed linking matrix represents the symmetric bilinear form~$Q$.
The~$4$-manifold~$X(L)$ formed by attaching $2$-handles along~$L$ (i.e. the \emph{trace on~$L$}) is a simply-connected~$4$-manifold with nonempty boundary and intersection form~$Q$.
Since~$Q$ is nonsingular,  the boundary of~$X(L)$ is an integral homology~$3$-sphere and therefore, by work of Freedman~\cite{Freedman}, bounds a contractible~$4$-manifold~$C$.
The closed simply-connected~$4$-manifold~$X:=X(L) \cup C$ now has intersection form~$Q$, as required.
When $Q$ is even (so that~$X$ is spin),  the Kirby-Siebenmann invariant is determined by the signature:~$\ks(X)=\sigma(X)/8$ mod~$2$; recall~\eqref{eq:KSSpin}.
%; this follows e.g. from the Kervaire-Milnor condition~\eqref{eq:KM} with $x=0$.
When $Q$ is odd (so that~$X$ is nonspin),  in order to realise both values of the Kirby-Siebenmann invariants,  one can use the star construction from~\cite[Section 10.4]{FreedmanQuinn} (see also~\cite{StongUniqueness} and~\cite[Section 5.11]{KasprowskiPowellRayCounterexamples}): this produces a~$4$-manifold homotopy equivalent to~$X$ but with the opposite Kirby-Siebenmann invariant.
\end{proof}

We conclude by emphasising the following point that was made at the end of the previous proof: closed simply-connected spin $4$-manifolds are homeomorphic if and only if their intersection forms are isometric.
Indeed, as was mentioned in~\eqref{eq:KSSpin},  for spin $4$-manifolds,  the Kirby-Siebenmann invariant is determined by the signature (and thus by the intersection form).

\subsection{A consequence of the sphere embedding theorem}

Both the $5$-dimensional~$s$-cobordism theorem for good groups and the fact that ``surgery theory works for good groups in dimension~$4$" rely on a deep result known as the disc embedding theorem; see~\cite{FreedmanQuinn,DET} for expositions.
Roughly speaking, the disc embedding theorem serves as a substitute for
the availability of the Whitney trick in high-dimensions.
% the high-dimensional Whitney trick
It allows for the high-dimensional proofs of the~$s$-cobordism theorem and of the surgery sequence to go through in dimension~$4$, when the fundamental group is good.
Another consequence of the disc (or sphere) embedding theorem is the following result whose underlying philosophy is similar: when the fundamental group is good, algebraic intersections can be realised be geometric intersections.
We refer to~\cite[Corollary~2.4]{PowellRayTeichner} for a proof.

\begin{theorem}
\label{thm:PRT}
Let~$V$ be a~$4$-manifold with good fundamental group~$\pi$.
If there is a hyperbolic form in~$(H_2(\widetilde{V}),\lambda_{\widetilde{V}})$, meaning that there is a~$\Z[\pi]$-submodule~$H\leq H_2(\widetilde{V})$ generated by a hyperbolic basis~$a_m,b_m$ for~$m=1,\ldots,k$ with
$$ \lambda_{\widetilde{V}}(a_i,b_j)=\delta_{ij} \quad \text{and} \quad \lambda_{\widetilde{V}}(a_i,a_j)=0=\lambda_{\widetilde{V}}(b_i,b_j) \quad \text{for every } i,j,$$
then there is a homeomorphism
$$V \cong W \#_{i=1}^k S^2 \times S^2.$$
\end{theorem}

With these theorems at hand, we now return to knotted surfaces.

\subsection{Applications to surfaces}

We describe how the~$s$-cobordism theorem,  surgery theory,  Freedman's classification of simply-connected~$4$-manifolds and Theorem~\ref{thm:PRT} can be leveraged to study knotted surfaces.
Throughout this section, we fix a closed simply-connected~$4$-manifold~$X$ and a group~$\pi$.
The take away from the following proposition is that uniqueness questions for knotted surfaces are related to the problem of extending a $3$-manifold homeomorphism to a homotopy equivalence on given bounding $4$-manifolds.

\begin{proposition}
\label{prop:scobApp}
Let~$S_0,S_1 \subset X$ be closed surfaces with knot group~$\pi$ and fix a homeomorphism of their closed tubular neighborhoods:
$$H \colon (\overline{\nu}(S_0), S_0) \to (\overline{\nu}(S_1),S_1).$$
%nd consider the homeomorphism~$h:=H| \colon \partial X_{S_0} \to X_{S_1}$.
If~$\pi$ is trivial,  finite cyclic, or a solvable Baumslag--Solitar group and~$H|_{\partial}$ extends to a homotopy equivalence~$X_{S_0} \to X_{S_1}$,  then the surfaces~$S_0$ and~$S_1$ are equivalent.
\end{proposition}
\begin{proof}
Apply Theorem~\ref{thm:SurgeryTheory} to the homotopy equivalence~$X_{S_0} \to X_{S_1}$ to deduce that~$X_{S_0}$ and~$X_{S_1}$ are~$s$-cobordant rel.~$H|_{\partial}$.
Since~$\pi$ is good, the relative~$s$-cobordism theorem implies that~$H|_{\partial}$ extends to a homeomorphism~$X_{S_0} \to X_{S_1}$.
Combining this homeomorphism with~$H$ leads to a homeomorphism~$(X,S_0) \to (X,S_1)$ thus establishing that the surfaces are equivalent.
\end{proof}

We move on to realisation.
 For this, we write~$N_{g,e}$ for the Euler number~$e$ disc bundle over the genus~$g$ surface and~$Y_{g,e}=\partial N_{g,e}$ for its boundary,  the Euler number~$e$ circle bundle over the genus~$g$ surface.
 In other words,~$Y_{g,e}$ is both the boundary of the closed tubular neighborhood of an Euler number~$e$ genus~$g$ surface~$S \subset X$ and the boundary of the exterior of this surface.
The take away from the following proposition is that realisation questions for knotted surfaces are related to the problem of building fillings of~$Y_{g,e}$.

\begin{proposition}
\label{prop:RealisationIdea}
Fix a group $\pi$.
If the circle bundle~$Y_{g,e}$ bounds a~$4$-manifold~$W$ that satisfies
\begin{enumerate}
\item the Kirby-Siebenmann condition:~$\ks(W)=\ks(X),$
\item the fundamental group condition:~$\pi_1(W) \cong \pi,$
\item the capping off condition:~$\widehat{W}:=W \cup N_{g,e}$ is simply-connected and satisfies~$Q_{\widehat{W}} \cong Q_X$, 
\end{enumerate}
then~$X$ contains a closed~$\pi$-surface~$S \subset X$.
\end{proposition}
\begin{addendum}
If, in addition, in the capping off condition, the isometry~$Q_{\widehat{W}} \cong Q_X$ is \emph{pointed},  i.e. if it maps the class of the core of~$N_{g,e}$ to a fixed homology class~$x \in H_2(X)$, then the closed orientable~$\pi$-surface~$S \subset X$ can be assumed to represent~$x$.
\end{addendum}
\begin{proof}
The additivity of the Kirby-Siebenmann invariant (see e.g.~\cite[Theorem 9.2]{FriedlNagelOrsonPowell}) and the Kirby-Siebenmann condition (1) imply that~$\ks(\widehat{W})=\ks(W)=\ks(X)$.
Since~$\widehat{W}$ is simply-connected (thanks to the capping off condition (3)) and has the same Kirby-Siebenmann invariant as~$X$, Freedman's classification of simply-connected~$4$-manifolds (Theorem~\ref{thm:HomeoClassification}) ensures that there is a homeomorphism~$\Phi \colon W \cong X$.
The required surface is then obtained by transporting the core of $N_{g,e} \subset \widehat{W}$ (i.e. the $0$-section) through $\Phi$.
The fundamental group condition (2) ensures that this surface has knot group $\pi$.
%as the~$0$-section of~$\Phi(\nu(S)) \subset X$.
For the addendum, one uses that the homeomorphism~$\Phi$ can be assumed to realise the given isometry of the intersection form~\cite{Freedman}.
%Wall probably gives a realisation by a smooth h.cob so don't need Quinn?
\end{proof}

Deciding whether or not there exists a filling $W$ of $Y_{g,e}$ with $\pi_1(W) \cong \pi$ is typically not overly challenging, as is deciding whether or not~$\widehat{W}$ can be arranged to be simply-connected.
On the other hand,  constructing $W$ so that~$Q_{\widehat{W}} \cong Q_X$ is usually more arduous.

We conclude with a first application of Theorem~\ref{thm:PRT} to the topic of representing a given homology class by a sphere with trivial knot group. 
\begin{proposition}
\label{prop:WarmUPLW}
If a homology class~$x \in H_2(X)$ is such that~$x \oplus 0 \in H_2(X\#_k S^2\times S^2)$ is represented by a sphere~$S \subset X\#_k S^2\times S^2$ with trivial knot group, then~$x$ is itself represented by a sphere with trivial knot group.
\end{proposition}
\begin{proof}
Set $X_k:=X\#_k S^2\times S^2$ and note the pointed isometry~$(Q_{X_k},[S]) \cong \left(Q_X \oplus \bsm 0&1 \\ 1&0\esm^{\oplus k},x \oplus 0\right)$.
Use~$X_S$ to denote the exterior of the sphere~$S \subset X_k.$
Since~$S$ represents~$x \oplus 0$,  a Mayer-Vietoris argument shows that the hyperbolic summand of~$H_2(X_k)$ determines a hyperbolic subgroup of~$H_2(X_S)$.
Theorem~\ref{thm:PRT} therefore implies that~$X_S=W \#_k S^2\times S^2$ for some~$4$-manifold~$W$.
This~$4$-manifold is a simply-connected filling of the circle bundle~$\partial X_S=Y_{e,0}$,  where~$e$ denotes the Euler number of~$S$.
A van Kampen argument shows that~$\widehat{W}=W \cup N_{e,0} $ must be simply-connected.
Thus, since~$\widehat{W}$ is simply-connected, has the same Kirby-Siebenmann invariant as~$X$,  and since there is a pointed isometry~$(Q_{\widehat{W}},[S])  \cong (Q_X,x)$ 
%Because~$(Q_{X_k},[S]) \cong (Q_X,x \oplus 0)$. but think through again.
(the addendum of) Proposition~\ref{prop:RealisationIdea} produces the required sphere in~$X$ representing~$x$.
\end{proof}

Proposition~\ref{prop:WarmUPLW} illustrates a fruitful idea in this area: by leveraging cancellation results in algebra or Theorem~\ref{thm:PRT},  there are situations where it is enough to prove a theorem stably, i.e.  up to taking sufficiently many connected sums  with $S^2 \times S^2$.
We refer to Lecture~\ref{sec:4}  and, more generally,  to~\cite{LeeWilczy,LeeWilczyOdd,LeeWilczyGenus,ConwayOrsonPencovitch,
ConwayGalvin,PencovitchSlice} for proofs involving this strategy.

 \section{Lecture three: Uniqueness results}
\label{sec:3}

Necessary conditions for two surfaces~$S_0,S_1 \subset X$ to be equivalent include that they share the same Euler number, that they both be ordinary/characteristic and that their exteriors have isometric equivariant intersection forms.
This latter condition can be refined,  namely by requiring that the isometry be compatible with a homeomorphism on the boundary of the surface exteriors.
We make this more precise, focusing on the case of $\Z_d$-surfaces for simplicity.
\begin{definition}
For $\Z_d$-surfaces~$S_0,S_1 \subset X$,  the \emph{boundary isomorphism} induced by an isometry~$F \colon \lambda_{\widetilde{X}_{S_0}} \cong \lambda_{\widetilde{X}_{S_1}}$ refers to the $\Z[\Z_d]$-isomorphism 
$$
\partial F \colon H_1(\partial \widetilde{X}_{S_0}) \xrightarrow{\cong} H_1(\partial \widetilde{X}_{S_1})
$$
that arises from the following commutative diagram:
$$
\xymatrix@C0.4cm@R0.4cm{
%H_2(\partial \widetilde{X}_{S_0})
%\ar[r]\ar[d]_{h_*}&
H_2(\widetilde{X}_{S_0})
\ar[r]\ar[d]_{F}&
{\overbrace{H_2(\widetilde{X}_{S_0},\partial \widetilde{X}_{S_0})}^{\cong H_2(\widetilde{X}_{S_0})^*}}
\ar[r]\ar[d]_{(F^*)^{-1}}&
H_1(\partial \widetilde{X}_{S_0})
\ar@{-->}[d]_{\partial F}\ar[r]&
0 \\
%%%%%%%%%%%%
%H_2(\partial \widetilde{X}_{S_1})
%%
%\ar[r]&
H_2(\widetilde{X}_{S_1})
\ar[r]&
{\underbrace{H_2(\widetilde{X}_{S_1},\partial \widetilde{X}_{S_1})}_{\cong H_2(\widetilde{X}_{S_1})^*}}
\ar[r]&
H_1(\partial \widetilde{X}_{S_1})\ar[r]&
0.
}
$$
Here, we write~$H_2(\widetilde{X}_{S_i})^*:=\Hom_{\Z[\Z_d]}(H_2(\widetilde{X}_{S_i}),\Z[\Z_d])$ for brevity and note that the left hand side square commutes because $F$ is an isometry.
\end{definition}

A suggested exercise involves thinking through the isomorphism $H_2(\widetilde{X}_{S_i},\partial \widetilde{X}_{S_i}) \cong H_2(\widetilde{X}_{S_i})^*$ (which is induced by Poincar\'e duality and the evaluation map) as well as showing that if the surfaces~$S_0$ and~$S_1$ are equivalent, then there is an isometry~$F \colon \lambda_{\widetilde{X}_{S_0}} \cong \lambda_{\widetilde{X}_{S_1}}$ and a bundle isomorphism $H \colon \overline{\nu}(S_0) \to \overline{\nu}(S_1)$ such that
$$\partial F=(H|_{\partial})_* \colon H_1(\partial \widetilde{X}_{S_0}) \to H_1(\partial \widetilde{X}_{S_1}).$$
The next theorem, which appears in~\cite[Theorem 1.5]{ConwayKasprowskiSurfaces},  shows that the converse nearly holds; throughout these lectures, we will treat $k$-invariants as a blackbox, referring to~\cite{ConwayKasprowskiKinvariant} for details.
We focus on closed orientable surfaces in closed~$4$-manifolds for simplicity, but note that the same result holds for nonorientable surfaces except when $X$ is spin and $d=1$; the statements can be further generalised to include surfaces with boundary; see also Remark~\ref{rem:UniquenessBoundary} below.

\begin{theorem}
\label{thm:CK}
Let~$X$ be a closed simply-connected~$4$-manifold and let~$S_0,S_1 \subset X$ be closed orientable~$\Z_d$-surfaces with the same Euler number that are either both characteristic or both ordinary.
%%Same genus follows from the conditions below.
The following assertions are equivalent:
\begin{itemize}
\item There is a homeomorphism~$(X,S_0) \to (X,S_1)$.
\item There is an isometry~$F\colon \lambda_{\widetilde{X}_{S_0}} \cong \lambda_{\widetilde{X}_{S_1}}$ and a bundle isomorphism~$H \colon \overline{\nu}(S_0) \to \overline{\nu}(S_1)$ such that~$(F,H|_{\partial})$ preserves the relative~$k$-invariants,
$$\partial F=(H|_{\partial})_* \colon H_1(\partial \widetilde{X}_{S_0}) \to H_1(\partial \widetilde{X}_{S_1}),$$
and such that when the~$S_i$ are characteristic,~$H|_{S_0}$ preserves the Rochlin quadratic forms.
\end{itemize}
The~$k$-invariant condition can be omitted when: either~$\pi$ is trivial,  or~$\pi$ is infinite cyclic,  or the~$S_i$ are spheres.
%, or the~$S_i$ are projective planes. 
The condition on the quadratic forms can be omitted for spheres.
\end{theorem}

We note that it is possible to reformulate Theorem~\ref{thm:CK} in terms of surface invariants so that it fits the mold of the uniqueness ``meta-statement" from Lecture~\ref{sec:1}; see~\cite[Theorem 1.9]{ConwayKasprowskiSurfaces}.
Next, we discuss  the philosophy underlying Theorem~\ref{thm:CK} and how it compares to previous results in the area.
Afterwards,  we outline some of the ideas of the proof.

\begin{remark}
\label{rem:Closed}
We compare Theorem~\ref{thm:CK} with prior literature on knotted surfaces and outline how it recovers most of the results from Lecture~\ref{sec:1}.
This is also an occasion to survey some of the results that have been proved in the area since the early successes that were listed during Lecture~\ref{sec:1}.
For brevity,  we refer to a pair $(F,H)$ satisfying~$\partial F=(H|_{\partial})_*$ as an \emph{extendable compatible pair.}
\begin{itemize}
\item When~$\pi$ is trivial and the surfaces are orientable,  Theorem~\ref{thm:CK} was (implicitly) proved by Boyer~\cite{BoyerRealization} with different methods.
In fact, starting from the outcome of Theorem~\ref{thm:CK},  Boyer establishes 
%that the intersection forms are isometric as well as 
the existence of an extendable compatible pair~$(F,H)$,  leading to the stronger result that homologous closed surfaces of the same genus with simply-connected complements are isotopic; recall Theorem~\ref{thm:Boyer}.
\item When~$\pi=\Z_d$ is finite cyclic and the surfaces are spheres, Theorem~\ref{thm:CK} can be shown to recover Lee and Wilczy\'nski's Theorem~\ref{thm:LW}; see~\cite[Section 7]{ConwayKasprowskiSurfaces}.
Here, the point is that a little work shows that given~$\Z_d$-spheres~$S_0,S_1 \subset X$, there exists an extendable compatible pair~$(F,H)$ if and only if the pointed hermitian forms~$(\lambda_{\Sigma_d(S_0)},[\widetilde{S}_0])$ and $(\lambda_{\Sigma_d(S_1)},[\widetilde{S}_1])$ are pointed isometric~\cite[Proposition 7.2]{ConwayKasprowskiSurfaces}.
%Given a nonzero divisibility~$d$ class~$x \in H_2(X),$ when~$X$ satisfies~$b_2(X)>6$ and
%\begin{equation}
%\label{eq:UniquenessIneq}
% b_2(X) > \underset{0 \leq j <d}{\operatorname{max}} \ \Big| \sigma(X)-\frac{2j(d-j)}{d^2}x\cdot x \Big|,
%\end{equation}
%Lee and Wilczy\'nski prove that such an isometry necessarily exists,  leading to the result that spheres representing~$x$ are isotopic.
%The condition~$b_2(X) >6$ can be replaced by one of the conditions in~\cite[Addendum 1]{LeeWilczy}.
%%When~$d>1$,  the inequality in~\eqref{eq:UniquenessIneq} precludes~$X$ from being definite.
\item When~$\pi \cong \Z$ is infinite cyclic,  Theorem~\ref{thm:CK} was proved in~\cite[Theorem 1.4]{ConwayPowell} with different methods.
Starting from the outcome of Theorem~\ref{thm:CK},  that article then used Blanchfield forms to show that an extendable compatible pair always exists,  so that closed~$\Z$-surfaces with isometric equivariant intersection forms are necessarily equivalent~\cite[Theorem~1.4]{ConwayPowell}.
Specifying to $X=S^4$ and combining with the discussion from Example~\ref{ex:Equivariant},  this leads to the result that any $\Z$-surface in~$S^4$ 
%of genus~$g \neq 1,2$ 
is unknotted,
%$\Z$-surfaces $S_0,S_1 \subset S^4$ of genus~$g \neq 1,2$ are isotopic
generalising Freedman's sphere-unknotting-Theorem~\ref{thm:FreedmanUnknotting}; see~\cite[Theorem~1.1]{ConwayPowell} and~\cite[Theorem 1.1]{ConwayUnknotting}.
\item When~$\pi=\Z_2,X=S^4$,  the surfaces are of nonorientable genus~$h$, and the Euler numbers agree and satisfy~$|e|\neq 2h$,  Theorem~\ref{thm:CK} follows from~\cite[Theorem A]{ConwayOrsonPowell}; the latter establishes that the surfaces are isotopic (i.e.  that most $\Z_2$-surfaces in $S^4$ are unknotted).
The same is true when~$h \leq 5$ and~$|e|=2h$ (with the case~$h=1$ originally due to Lawson~\cite{Lawson} and the cases~$h=4,5$ requiring~\cite{Pencovitch}).
We also refer to~\cite{FinashinKreckViro,KreckOnTheHomeomorphism} for other early work in the area.
A suggested exercise involves showing that Theorem~\ref{thm:CK} recovers Lawson's result (stated in Theorem~\ref{thm:LawsonUnknotting}) that~$\Z_2$-projective planes in~$S^4$ are unknotted.
We also refer to Example~\ref{ex:Z2SpheresCP2} for an application to~$\Z_2$-spheres in $\C P^2$ (i.e.  for how Theorem~\ref{thm:CK} relates to~\cite{ConwayOrson}) and to~\cite[Theorem 1.1]{ConwayKasprowskiSurfaces} for %further applications, e.g.
another application, namely
 for the fact that any knot in $S^3$ whose determinant is a prime power bounds at most one~$\Z_2$-Moebius band in the~$4$-ball.
%%Don't delete
%%If can get ride of k-invariant, then can write something like the following:
%When the~$e=2h$ with~$h>6$, our results are new.
\end{itemize}
In summary,  while Theorem~\ref{thm:CK} might not provide the sharpest statement in any particular case, it unifies all previous results in the area.
It should therefore be viewed as providing a common framework from which further uniqueness statements can then be deduced.
%In addition,  the previous discussion shows that Theorem~\ref{thm:CK} recovers
\end{remark}

We now outline the proof of Theorem~\ref{thm:CK}.
\begin{proof}[Proof outline of Theorem~\ref{thm:CK}]
As noted above, if the surfaces are equivalent, then an extendable compatible pair exists (and some additional work establishes the $k$-invariant and spin conditions); we therefore focus on proving the converse.
Fix a bundle isomorphism~$H \colon \nu(S_0) \to \nu(S_1)$ and an isometry $\lambda_{\widetilde{X}_{S_0}} \cong \lambda_{\widetilde{X}_{S_1}}$ as in the theorem.
This first step of the proof is to adapt recent work on 4-manifolds with boundary~\cite{ConwayKasprowski4Manifolds} to show that $H|_\partial$ extends to a homotopy equivalence if the three following conditions hold: (a) the pair~$(F,H|_{\partial})$ preserves the relative~$k$-invariants,  (b) the pair~$(F,H|_{\partial})$ fits into the commutative diagram
\begin{equation}
\label{eq:CompatiblePair}
\xymatrix@C0.4cm@R0.4cm{
H_2(\partial \widetilde{X}_{S_0})
\ar[r]\ar[d]_{(H|_{\partial})_*}&
H_2(\widetilde{X}_{S_0})
\ar[r]\ar[d]_{F}&
{\overbrace{H_2(\widetilde{X}_{S_0},\partial \widetilde{X}_{S_0})}^{\cong H_2(\widetilde{X}_{S_0})^*}}
\ar[r]\ar[d]_{(F^*)^{-1}}&
H_1(\partial \widetilde{X}_{S_0})
\ar[d]_{(H|_{\partial})_*}\ar[r]&
0 \\
%%%%%%%%%%%%
H_2(\partial \widetilde{X}_{S_1})
\ar[r]&
H_2(\widetilde{X}_{S_1})
\ar[r]&
{\underbrace{H_2(\widetilde{X}_{S_1},\partial \widetilde{X}_{S_1})}_{\cong H_2(\widetilde{X}_{S_1})^*}}
\ar[r]&
H_1(\partial \widetilde{X}_{S_1})\ar[r]&
0,
}
\end{equation}
and (c),  when the universal covers of $X_{S_0}$ and $X_{S_1} $ are spin,  the universal cover of~$X_{S_0} \cup_{H|_\partial} X_{S_1}$ is also spin.
The second step of the proof involves investigating when (a) and (b) hold and showing that~$(F,H)$ can be modified if necessary so as to yield a new 
%extendable compatible 
pair $(F',H')$ for which the universal cover of the union is spin when relevant (i.e.  arranging condition (c)) without affecting~(a) and~(b).
When $(a),(b)$ and $(c)$ can be achieved, combining these two steps ensures that $H|_\partial$ extends to a homotopy equivalence $X_{S_0} \to X_{S_1}$,  and Proposition~\ref{prop:scobApp} then implies that the surfaces~$S_0$ and~$S_1$ are equivalent.

We focus on outlining (c) of the second step; indeed imposing (a) for spheres involves technical results on the relative~$k$-invariant whose definition we have not introduced, whereas establishing~(b), i.e. showing that the left hand square in~\eqref{eq:CompatiblePair} commutes,  relies on somewhat routine arguments from algebraic topology.
For (c), firstly, note that the exterior~$X_S$ of a surface~$S \subset X$ is spin if and only if either~$X$ is spin or~$S$ is characteristic; see e.g.~\cite[Proposition 7.18]{FriedlNagelOrsonPowell}.
Secondly,  if~$X_S$ is nonspin but~$\widetilde{X}_S$ is spin, then it can be shown that~$d$ is necessarily even.
There are therefore three cases to consider.
When~$X$ is spin,  the union~$X_{S_0} \cup_{H|_{\partial}} X_{S_1}$ may not be spin but we show it is possible to change~$(F,H)$ in order to obtain a new pair~$(F',H')$ such that~$(F',H'|)$ still preserves~$k$-invariants, satisfies~$\partial F'=(H'|_\partial)_*$ and such that, additionally,  the union~$X_{S_0} \cup_{H'|_{\partial}} X_{S_1}$ is spin.
When~$d$ is even, we show that the universal cover of the union is spin for any~$H$.
When~$S$ is characteristic,  we show that~$X_{S_0} \cup_{H|_{\partial}} X_{S_1}$ is spin if and only if~$H|_{S_0}$ preserves the quadratic~forms.

As explained above,  thanks to this outline, $H|_\partial$ extends to a homotopy equivalence $X_{S_0} \to X_{S_1}$,  and Proposition~\ref{prop:scobApp} then implies that the surfaces~$S_0$ and~$S_1$ are equivalent.
\end{proof}

Next, we apply Theorem~\ref{thm:CK} to~$\Z_2$-spheres in $\C P^2$.
We refer to the exercises for how Theorem~\ref{thm:CK} recovers Lawson's result that $\Z_2$-projective planes in $S^4$ are unknotted.

\begin{example}
\label{ex:Z2SpheresCP2}
We use Theorem~\ref{thm:CK} to provide a quick proof that~$\Z_2$-spheres in~$\C P^2$ are equivalent;
% up to equivalence, there is a unique~$\Z_2$-sphere in~$\C P^2$; 
 this recovers most of~\cite[Theorem 1.1]{ConwayOrson} with the caveat that the latter article then upgrades the result to an isotopy classification.
Given two $\Z_2$-spheres~$S_0,S_1 \subset \C P^2$ (which are necessarily ordinary and have Euler number~$4$),  in order to prove equivalence, by Theorem~\ref{thm:CK}, it suffices to find an isometry~$F \colon \lambda_{\widetilde{X}_{S_0}} \cong \lambda_{\widetilde{X}_{S_1}}$ and a bundle isomorphism~$H \colon \nu(S_0) \to \nu(S_1)$  such that~$\partial F=(H|_{\partial})_*$ on~$H_1(\partial \widetilde{X}_{S_0})$.
In fact,  since~$H_1(\partial \widetilde{X}_{S_0}) \cong \Z_2$ admits a unique automorphism,  this latter condition is automatic.
% by flipping the sign of~$F$ if necessary.
But now, the existence of an isometry~$\lambda_{\widetilde{X}_{S_0}} \cong \lambda_{\widetilde{X}_{S_1}}$ follows from Example~\ref{ex:Equivariant}  which asserted that the equivariant intersection form of any~$\Z_2$-sphere~$S\subset \C P^2$ is isometric to~$(\Z_-,(x,y) \mapsto 2(1-T)xy)$.
This shows that,  up to equivalence,
%and orientation reversal, 
there exists a unique $\Z_2$-sphere in $\C P^2$,  namely the algebraic curve given by~$x^2 + y^2+z^2=0$.
In fact, with a little additional work,  it is possible to show that homologous $\Z_2$-spheres in $\C P^2$ are isotopic~\cite[Theorem 1.1]{ConwayOrson}.

Combined with the work of Lee and Wilczynski~\cite{LeeWilczy} and~\cite[Theorem 1.4]{ConwayPowell},  it then follows that a class~$x \in H_2(\C P^2) \cong \Z$  is represented by
a locally flat sphere with abelian knot group if and only if~$|d|=0,1,2$ and that this sphere is unique up to isotopy~\cite[Theorem 1.2]{ConwayOrson}.
In other words, up to isotopy, there are precisely~$5$ spheres in~$\C P^2$ with abelian knot group.
\end{example}

We conclude this section by surveying 
%the state of the art concerning 
uniqueness results for surfaces with boundary.
\begin{remark}
\label{rem:UniquenessBoundary}
Several of the articles listed in Remark~\ref{rem:Closed} also allow for properly embedded surfaces with boundary a knot in $S^3$.
For uniqueness results on~$\Z$-surfaces,  we refer to~\cite{ConwayPowell}; see also~\cite{ConwayDaiMiller} for the fact that the equivariant intersection form is not a complete invariant when the boundary is nonempty.
For surfaces with boundary and simply-connected complement, we refer to~\cite{Pyronneau}; see also~\cite{PencovitchSlice}.
For further results on discs in the~$4$-ball, see~\cite{ConwayPowellDiscs,ConwayDiscs,ConwayUnknotting}, whereas we refer to~\cite{ConwayOrsonPowell} for $\Z_2$-surfaces in the $4$-ball.
Again,~\cite{ConwayKasprowskiSurfaces} provides a framework unifying these results.
We also mention~\cite{OrsonPowellSpines} for spheres in knot traces and~\cite{ConwayMiller} for immersed~$\Z$-surfaces.
\end{remark}

\section{Lecture four: existence results}
\label{sec:4}

In order to establish existence results, the idea is nearly always to ``build the surface exterior" (i.e. to build an appropriate filling for the circle bundle~$Y_{e,g}$) as mentioned at the end of Lecture~\ref{sec:2}; see e.g.~\cite{KasprowskiPowellRayTeichnerEmbedding} for a different approach.
Contrarily to uniqueness results however, there is currently no single approach that encompasses all previous results in the area.
Broadly speaking,  here are the two strategies that have been used to build closed surfaces in closed $4$-manifolds: 
\begin{itemize}
\item The first is to attempt to \emph{stably} build  the surface (or the filling) and to then use Theorem~\ref{thm:PRT} to obtain the required surface~\cite{LeeWilczyOdd,LeeWilczy, LeeWilczyGenus,ConwayGalvin}. 
All of these papers involve surfaces with finite cyclic knot groups.
\item The second strategy involves building a cobordism from~$Y_{e,g}$ to a~$3$-manifold with very simple algebraic topology which can then be capped off to obtain the required filling; see~\cite{BoyerRealization} for surfaces with trivial knot group and~\cite{ConwayPiccirilloPowell} for surfaces with knot group~$\Z$.
\end{itemize}

We will focus on the first strategy since it also encompasses realisation results for surfaces with trivial knot group. 
In fact, at the time of writing, the first strategy recovers all known realisation results besides those involving $\Z$-surfaces.
For an expository account of the second strategy in the context~$\Z$-surfaces, we refer to~\cite{ConwayOberwolfach,ConwayMontreal}.

\begin{remark}
For surfaces with boundary,  the first strategy was used in~\cite{ConwayOrsonPencovitch,PencovitchSlice}, whereas the second featured in~\cite{ConwayPiccirilloPowell}; see also~\cite{ConwayMiller}.
There is also a third strategy which consists of first making the boundary of the putative surface exterior bound a (Poincar\'e) complex with the desired homotopy type and then using surgery theory to upgrade this filling to a (manifold) filling~\cite{Freedman,FreedmanQuinn,FriedlTeichner,Davis,FellerMillerNagelOrsonPowellRay}.
At the time of writing,  this strategy has only been successfully applied in settings where the surface exteriors are aspherical.
\end{remark}

Next, we outline the proof of Lee and Wilczy\'nski's Theorem~\ref{thm:LW} in the case of spheres.
The idea for higher genus surfaces is similar but technically more demanding.

\begin{theorem}
\label{thm:LWRealisation}
Given an integer~$d>0$, a class~$x \in H_2(X)$ is represented by a $\Z_d$-sphere if and only if the divisibility of $x$ equals $d$,
\begin{equation}
\label{eq:ExistenceLW}
 b_2(X) \geq  \underset{0 \leq j <d}{\operatorname{max}} \ \Big| \sigma(X)-\frac{2j(d-j)}{d^2}x\cdot x \Big| \tag{Ro}
 %%Note that here it's an inequality.
 \end{equation}
and, when~$x$ is characteristic,  it additionally satisfies
\begin{equation}
\label{eq:ExistenceCharacteristic}
\ks(X)=\frac{1}{8}(\sigma(X)-x\cdot x) \mod 2. \tag{km}
\end{equation}
\end{theorem}
\begin{proof}[Proof outline of Theorem~\ref{thm:LWRealisation}]
The necessity of the conditions can be verified using work of Rochlin~\cite{Rochlin} for~\eqref{eq:ExistenceLW} and Freedman-Kirby~\cite{FreedmanKirby} for~\eqref{eq:ExistenceCharacteristic}.
We therefore focus on sufficiency.
The key idea is contained in the following claim which generalises Proposition~\ref{prop:WarmUPLW}.
\begin{claim}
If there exists a~$\Z_d$-sphere~$S \subset X \#_k S^2 \times S^2$ for some~$k \geq 0$ such that
\begin{itemize}
\item the sphere represents the class~$x \oplus 0$ and 
\item there is a pointed isometry
$$(\lambda_{\Sigma_d(S)},[\widetilde{S}]) \cong (\lambda',z') \oplus \left(\bsm 0&1 \\ 1&0 \esm^{\oplus k},0\right),$$
for some pointed hermitian form $(\lambda',z')$ such that~$(\lambda',z') \otimes_{\Z[\Z_d]}\Z \cong (Q_X,x),$
\end{itemize}
then~$x$ is realised unstably by a sphere, i.e. the conclusion of the theorem holds.
 \end{claim}
 \begin{proof}[Proof of the claim]
 Informally, the goal is to surger away the hyperbolic summand without damaging the sphere~$S$.
Writing~$\widetilde{X}_S$ for the universal cover of~$X_S$ so that~$\Sigma_d(S)=\widetilde{X}_S \cup \overline{\nu}(S)$,  a Mayer-Vietoris argument shows that the~$0 \oplus \Z[\Z_d]^{2k}$-summand lives in~$H_2(\widetilde{X}_S) \cong \pi_2(\widetilde{X}_S) \cong \pi_2(X_S)$.
Since~$\pi_1(X_S)\cong \Z_d$ is a good group,  Theorem~\ref{thm:PRT} then ensures  that $X_S \cong W \#_k S^2 \times S^2$ for some~$4$-manifold $W$.
By construction, this latter~$4$-manifold has~$\pi_1(W)\cong \Z_d$ and boundary the circle bundle~$Y_{e,0}$.
We now apply the criterion from (the addendum to) Proposition~\ref{prop:RealisationIdea}.
Recall that~$N_{e,0}$ denotes the~$D^2$-bundle over the (abstract) sphere~$\Sigma$, so that the closed~$4$-manifold~$\widehat{W}:=W \cup N_{e,0}$ is simply-connected,  satisfies~$\ks(\widehat{W})=\ks(W)=\ks(X_S)=\ks(X)$,  and there is a pointed isometry
%%The branched cover along \Sigma has the hyperbolics chopped off so the intersection form is \lambda'.
$$(Q_{\widehat{W}},[\Sigma])
 \cong (\lambda_{\Sigma_d(\Sigma)},[\widetilde{\Sigma}])  \otimes_{\Z[\Z_d]}\Z
 \cong (\lambda',z')  \otimes_{\Z[\Z_d]}\Z
 \cong (Q_X,x)$$ 
Proposition~\ref{prop:RealisationIdea} now implies that there is a sphere in $X$ representing $x$.
 \end{proof}

It remains to outline how to arrange for the conditions of the claim to hold.
The idea of the proof is to first  show that the condition on the Arf invariant guarantees that,  for some~$k \geq 0$, the class~$x \oplus 0$ is represented by a $\Z_d$-sphere~$S \subset X \#_k S^2 \times S^2$~\cite[Theorem~2.1]{LeeWilczyOdd}, and to then use the main hypothesis of the theorem (namely~\eqref{eq:ExistenceLW}) to establish the splitting condition.

We describe the idea underlying the first step as the second involves arduous algebra; we refer to~\cite{LeeWilczyOdd,LeeWilczy} for the details.
Represent~$x \in H_2(X)$ by a surface~$S \subset X$.
The goal is to show that,  assuming~\eqref{eq:ExistenceCharacteristic} in the characteristic case,  it is possible to ambiently surger~$S$ to a sphere at the price of stabilising~$X$.
Consider the task of performing ambient surgery on an embedded loop~$\gamma \subset S$.
This requires the existence of an embedded disc~$D \subset X$ with~$\partial D=\gamma=D \cap S$ and such that disc framing on~$\gamma$ agrees with the surface-induced framing on it; see e.g.~\cite[page~508]{Scorpan} for details and illustrations.
%Cite Freedman-Kirby instead of FQ for boundary twist.
Write~$W(D)\in \Z$ for the obstruction to this framing condition being realised. 
It is not too challenging to show that if~$W(D)$ and~$D \cdot S$ have the same parity, then ambient surgery along~$\gamma$ is possible in a stabilisation of~$X$.
%\begin{claim}
%If~$W(D)$ and~$D \cdot S$ have the same parity, then ambient surgery along~$\gamma$ is possible in a stabilisation of~$X$.
%\end{claim}
%\begin{proof}
%Performing a so-called \emph{boundary twist} on~$D$ (see e.g.~\cite[p.~15]{FreedmanQuinn}) changes both quantities by~$\pm 1$ so that, after sufficiently many such operations, it is possible to assume that~$D \cdot S=0$ and that~$W(D)$ is even.
%It is then possible to change~$W(D)$ by any even number (without affecting~$D \cdot S)$ at the price of stabilisations; see e.g.~\cite[page 93]{FreedmanKirby} or~\cite[page 150]{Scorpan}.
%Thus we can assume that both~$W(D)$ and~$D \cdot F$ vanish.
%The Whitney trick can be performed stably (see e.g.~\cite[page 150]{Scorpan}), thus making it possible to arrange that~$D \cap F=\gamma$, as needed.
%% to prove the claim.
%\end{proof}

It therefore remains to find an embedded disc~$D$ such that~$W(D)+D\cdot S$ is even.
It turns out that if the class~$x$ is ordinary, this is always possible, whereas in the characteristic case, the quantity~$q_S([\gamma]):=W(D)+D\cdot S$ leads to a genuine obstruction.
Thus, in the ordinary case,  it is possible to iteratively stably ambiently surger the curves in a symplectic basis of~$H_1(S)$ leading to a sphere~$S$, whereas it turns out that in the characteristic case, this is also possible provided~$\operatorname{Arf}(q_S)=0$, i.e. (using Remark~\ref{rem:Rochlin}) provided~$\ks(X) -1/8(\sigma(X)-x \cdot x)=0$. 
This concludes the proof outline that $x \oplus 0$ is represented by a sphere $S'$.

Since $x$ has divisibility $d$,  the sphere $S'$ satisfies $H_1(X_{S'})\cong \Z_d$ but  it might not have knot group $\Z_d$.
To correct this,  one can perform additional stabilisations:
%to arrange for~$x \oplus 0$ to be stably representable by a~$\Z_d$-sphere
surger push offs of curves representing the commutators of the fundamental group of the sphere exterior.
This concludes our outline of the first step of the proof and thus our outline of Theorem~\ref{thm:LWRealisation}.
\end{proof}

Lee and Wilczy\'nski's results show how to realise homology classes by surfaces, but there are far fewer results whose focus is on realising further invariants such as the equivariant intersection form.
To the author's knowledge,  this only occurs in~\cite{ConwayPiccirilloPowell} and~\cite{ConwayGalvin}.
%Boyer also maybe.
We outline the proof of this latter article since it also uses the aforementioned ``stable strategy" and leads both to first example of a nonsmoothable surface in $S^4$ and to the first example of a nontrivial knot in $S^4$ with cyclic knot group, thus answering~\cite[Questions 4.29(a) and 4.32]{K3}.

\begin{theorem}
\label{thm:ConwayGalvin}
For every $n\geq 1$, there is a~$\Z_2$-surface~$S \subset S^4$ whose double branched cover has intersection form
$$Q_{\Sigma_2(S)} \cong E_8 \oplus (1)^{\oplus n}.$$
In particular, there exists infinitely many nonsmoothable surfaces in the~$4$-sphere.
\end{theorem}

Thus, Theorem~\ref{thm:ConwayGalvin} disproves the ``nonorientable unknotting conjecture" that every $\Z_2$-surface in~$S^4$ is unknotted.
More precisely,  using $(e,h)$ to denote the pair consisting of the Euler number and the nonorientable genus, the nonorientable unknotting conjecture holds when $h \leq 5$ or when~$|e| \neq 2h$ (recall Remark~\ref{rem:Closed}) but is false for every~$h \geq 9$ (and~$|e|=2h$).
%The cases where~$h=6,7,8$ and~$|e|=2h$ remain open.
%In the orientable case,  as mentioned in Remark~\ref{rem:Closed}, the orientable unknotting conjecture is settled whenever the genus of the surface satisfies $g \neq 1,2$.
%that that every $\Z$-surface in~$S^4$ is unknottedis settle for every $g \neq 

\begin{proof}[Proof outline of Theorem~\ref{thm:ConwayGalvin}]
Firstly, we note that if the double branched cover of a knotted surface~$S \subset S^4$ has (nondiagonalisable) intersection form~$Q_{\Sigma_2(S)} \cong E_8 \oplus (1)^{\oplus n}$, then Donaldson's theorem implies that~$S$ is nonsmoothable~\cite{Donaldson}.
Secondly, if~$S \subset S^4$ satisfies~$Q_{\Sigma_2(S)} \cong E_8 \oplus (1)$ and~$U_k$ denotes the unknotted surface of nonorientable genus~$k$ and Euler number~$-2k$, then the new knotted surface~$S':=S \# U_k$ satisfies~$Q_{\Sigma_2(S')} \cong E_8 \oplus (1)^{\oplus (k+1)}$.

We therefore focus on outlining the construction of a surface~$S$ with~$Q_{\Sigma_2(S)} \cong E_8 \oplus (1)$.
Given a~$\Z_2$-surface~$S \subset S^4$ of nonorientable genus~$h$, it is possible to show that the equivariant intersection form~$\lambda_{\widetilde{X}_S}$ can be read off from the intersection form of the branched cover~$Q_{\Sigma_2(S)}$ and that the converse is true when~$h$ is not congruent to~$4$ mod~$8$.
In other words, if a symmetric bilinear form~$b$ is realisable by a~$\Z_2$-surface, i.e. if~$b \cong Q_{\Sigma_2(S)}$ (which implies that $b$ is nonsingular and odd),  then an associated~$\Z[\Z_2]$-hermitian form~$\lambda_b$ is also realised as~$\lambda_b\cong \lambda_{\widetilde{X}_S}$ and the converse holds in most cases.
Our goal is therefore to show that~$\lambda_b$ is realisable.
In fact, using~$Y_{e,h}$ to denote the Euler number~$e=-2h$ circle bundle over the nonorientable genus~$h$ surface~$\Sigma$,  we have the following.

\begin{claim}
\label{claim:StablyRealisable}
Let $b$ be a nonsingular odd symmetric bilinear form of rank $h \not\equiv 4$ mod $8$.
If~$Y_{e,h}$ bounds a~$4$-manifold $X$ with~$\pi_1(X) \cong \Z_2$,  such that the $S^1$-fibre of $Y_{e,h}$ generates~$\pi_1(X) \cong \Z_2$ and so that the equivariant intersection form of $X$ is stably isometric to~$\lambda_b$, then~$b$ is realisable by a~$\Z_2$-surface~$S \subset S^4$.
\end{claim}
\begin{proof}[Proof of the claim]
Assume there is a filling~$V$ of~$Y_{e,h}$ with~$\pi_1(V) \cong \Z_2$,  and an isometry
$$(H,\lambda_b)   \oplus \left( \Z[\Z_2]^2, \bsm 0&1\\ 1&0 \esm\right)^{\oplus k} \cong (H_2(\widetilde{V}),\lambda_{\widetilde{V}}).$$
Since finite groups are good,
%(see e.g.~\cite[Theorem 19.2]{DET}) 
%%Said so before
Theorem~\ref{thm:PRT} implies that~$V$ splits off~$k$ copies of~$S^2 \times S^2$.
Surgering these connected summands yields a filling $W$ of~$Y_{e,h}$ that realises~$(H,\lambda_b)$.

This way,~$\widehat{W}:=W \cup N_{e,h}$ (where~$N_{e,h}$ denotes the Euler number~$e=-2h$ disc bundle over the nonorientable genus~$h$ surface) is a closed simply-connected~$4$-manifold that contains a nonorientable genus~$h$ surface.
We now wish to apply Proposition~\ref{prop:RealisationIdea} to conclude.
Since we are trying to prove that $\widehat{W}$ is homeomorphic to~$S^4$,  it suffices to show that~$b_2(\widehat{W})=0$ (instead of computing the intersection form and Kirby-Siebenmann invariant) and this follows from a quick Euler characteristic calculation.
Thus~$\widehat{W}$ is homeomorphic to~$S^4$ (by Freedman's work~\cite{Freedman}) and contains a~$\Z_2$-surface~$S$ with non-orientable genus~$h$,  Euler number~$-2h$, and
exterior having equivariant intersection form isometric to~$\lambda_b$. 
As alluded to above,  since~$h \neq 8$ mod~$4$, it follows that~$Q_{\Sigma_2(S)}\cong b$.
This concludes the proof of the claim.
\end{proof}

It remains to outline how to arrange for the conditions of the claim to hold, i.e. for~$\lambda_b$ to be stably realisable.
Any odd positive-definite nonsingular form~$b$ is stably isometric to~$Q_{\Sigma_2(U)}=(1)^{\oplus h}$, the intersection form of the branched cover of the unknotted surface~$U \subset S^4$ of nonorientable genus~$h$ and Euler number~$-2h$.
When~$b=E_8 \oplus (1)$,  using a description of~$\lambda_b$ and~$\lambda_{\widetilde{X}_U}$ as pullback forms, it is then possible to construct a stable isometry~$\lambda_b \cong_s \lambda_{\widetilde{X}_U}$; this is the technical core of the argument, we refer to~\cite[Section 5]{ConwayGalvin} for the details.
This shows that~$\lambda_b$ is stably realisable, therefore ensuring that the conditions of the claim hold.
The claim then produces the desired~$\Z_2$-surface realising $b$.
\end{proof}

\color{black}

\section*{Exercises}

\begin{enumerate}
\item Describe the boundary of the exterior of a closed surface in a closed~$4$-manifold.
\item Calculate the homology of sphere exteriors in a closed, connected, simply-connected~$4$-manifold~$X$; deduce that (a) if~$\pi_1(X \setminus S)$ is abelian, then it is cyclic and (b) if~$S$ additionally has divisibility~$d$, then~$\pi_1(X \setminus S) \cong \Z_d$.
%\item Use Rochlin ineq and KM to study the genus function, then use the LW result.
\item For~$S \subset \C P^2$ a~$\Z_2$-sphere, calculate the intersection form of the universal cover of the exterior without drawing handle diagrams.
Along the way,  show that~$H_1(\partial \widetilde{X}_S)\cong \Z_2$.
\item Same question for~$S \subset S^4$ a~$\Z_2$-projective plane.
Use~$H_1(\partial \widetilde{X}_S)=\Z_4$ without proof.
\item Prove that if~$X_0$ and $X_1$ are $s$-cobordant rel.~$f$ and the relative $s$-cobordism theorem holds, then $f \colon \partial X_0 \to \partial X_1$ extends to a homeomorphism $X_0 \to X_1$.
%\item Adapt Propositions~\ref{prop:scobApp} and~\ref{prop:RealisationIdea} to surfaces with boundary a fixed knot.
\item Describe the isomorphism $H_2(\widetilde{X}_S,\partial \widetilde{X}_S) \cong \Hom_{\Z[\Z_d]}(H_2(\widetilde{X}_{S_i}),\Z[\Z_d])$ that was mentioned above Theorem~\ref{thm:CK}.
\item Prove that if two~$\Z_d$-surfaces $S_0,S_1 \subset X$ are equivalent, then there exists an isometry~$F \colon \lambda_{\widetilde{X}_{S_0}} \cong \lambda_{\widetilde{X}_{S_1}}$ and a bundle isomorphism $H \colon \overline{\nu}(S_0) \to \overline{\nu}(S_1)$ such that
$$\partial F=(H|_{\partial})_* \colon H_1(\partial \widetilde{X}_{S_0}) \to H_1(\partial \widetilde{X}_{S_1}).$$
\item Use Theorem~\ref{thm:CK} to prove Lawson's theorem that~$\Z_2$-projective planes in~$S^4$ are unknotted.
You can use without proof that for a $\Z_2$-surface $S \subset S^4$, the group~$\Z_2=\langle T \mid T^2=1\rangle$ acts on $H_2(\widetilde{X}_S)$ by $Tx=-x$ for every $x \in H_2(\widetilde{X}_S)$.
\end{enumerate}

\def\MR#1{}
\bibliography{BiblioPCMI}
\end{document}